\documentclass[12pt]{amsart}

\title[On nef reductions]{On nef reductions of  projective irreducible
       symplectic manifolds}
\author{Daisuke Matsushita}
\address{Division of Mathematics, Graduate School of Science,
         Hokkaido University,  Sapporo, 060-0810 Japan}
\thanks{* Partially supported by Grand-in-Aid \# 15740002
 (Japan Society for Promortion of Sciences).} 
\email{matusita@math.sci.hokudai.ac.jp}

\usepackage{amsfonts}
\usepackage{amsmath}
\usepackage{latexsym}
\usepackage[all]{xy}
\usepackage[dvipdfm]{hyperref}

\newtheorem{thm}{{\sc Theorem}}[section]

\newtheorem{defn}[thm]{{\sc Definition}}
\newtheorem{claim}[thm]{{\sc Claim}}
\newtheorem{prob}[thm]{{\sc Problem}}
\newtheorem{defnprop}[thm]{{\sc Definition--Theorem}}
\newtheorem{remark}[thm]{{\sc Remark}}
\theoremstyle{definition}
\newtheorem{say}[thm]{}

\begin{document}
\maketitle

\begin{abstract}
  Let $X$ be a projective irreducible symplectic manifold
  and $L$ is a non trivial nef divisor on $X$.
  Assume that the nef dimension of $L$
  is strictly less than the dimension of $X$.
  We prove that $L$ is semiample.
\end{abstract}

\section{Introduction}
 We begin with the definitions of {\itshape irreducible symplectic manifolds}
 and
 {\it  Lagrangian fibrations\/}.
\begin{defn}\label{defn_symplectic_variety}
  A  compact K\"{a}hler manifold $X$ is said to be
  a symplectic manifold if $X$ carries
  a holomorphic symplectic form.
  Moreover if   $X$ satisfies the following 
  two conditions, $X$ is said to be irreducible.
\begin{enumerate}
 \item $\pi_1 (X) = \{1\}$.
 \item $h^2 (X, \mathcal{O}_X) = 1$.
\end{enumerate}
\end{defn}

\begin{defn}\label{lagrangian_fibration}
 Let $X$ be a symplectic manifold,
 $\omega$ a symplectic form on $X$
 and $S$ a normal variety. A proper surjective morphism 
 with connected fibres
 $f \colon (X,\omega ) \to S$ 
 is said to be a Lagrangian fibration
 if a general fibre $F$ of $f$ is a Lagrangian variety with
 respect to $\omega$,
 that is,  $\dim F = (1/2)\dim X$
 and the restriction of the symplectic 2-form 
 $\omega |_{F}$ is 
 identically zero.
\end{defn}

 The simplest example of an irreducible symplectic manifold
 is  a $K3$ surface. It is expected that
 $K3$ surfaces and irreducible symplectic
 manifolds share many geometric properties.
 Let $L$ be a non trivial nef divisor on a $K3$ surface.
 It is well known that if $L^2 = 0$,
 then $L$ defines an elliptic fibration, which
 is the simplest example of a Lagrangian fibration.
 By using a Beauville-Bogomolov-Fujiki form $q$
 \cite[Th\'eor\`eme 5]{beauville},
 we can consider the following problem
 for an irreducible symplectic manifold,
 which is posed by Huybrechts and Sawon:

\begin{prob}
 Let $X$ be an irreducible symplectic manifold
 and $L$  a non trivial nef divisor on $X$.
 If $q(L) = 0$, then $L$ 
 defines a Lagrangian fibration?
\end{prob}

  We give a partial answer for the above problem. To state our
 results, we introduce the nef dimension of a nef divisor,
 which is due to \cite{eckl}.

\begin{defnprop}\label{Nef_Dimension}
\cite[Theorem 2.1 and Definition 2.7]{eckl}

 Let $X$ be a normal projective variety and $L$ a nef divisor on $X$.
 Then there exists a rational map $f : X \dasharrow S$
 which satisfies the following three conditions:
\begin{enumerate}
 \item A general fibre of $f$ is compact.
 \item $L$ is numerically trivial on a general fibre of $f$.
 \item For every general point $x \in X$ and every
       irreducible curve passing through $x$
       with $\dim f(C) > 0$, we have $L.C > 0$.
\end{enumerate}
 Moreover the map $f$ is unique up to birational equivalence of $S$.
 We define the nef dimension of $L$ as
$$
  n(L) := \dim S.
$$
\end{defnprop}

 Our result is the following:
\begin{thm}\label{main}
 Let $X$ be a projective irreducible symplectic manifold and
 $L$ a non trivial nef divisor on $X$.
 Assume that $n(L) < \dim X$. Then $L$ is semiample,
 that is, there exists an integer $M$ such that
 the linear system $|ML|$ is free.
\end{thm}

\begin{remark}
 Since the morphism $f : X \to S$ which is defined by
 $|ML|$ satisfies the three conditions of Definition--Theorem
 \ref{Nef_Dimension}, we have $\dim S = n(L)$.
 By \cite[Theorem 2]{matsu} and \cite[Theorem 1]{matsu},
 $f$ is a Lagrangian fibration.
\end{remark}

\begin{remark}\label{relation_of_dimensions}
 Let $L$ be a nef divisor on a projective variety.
 By \cite[Proposition 2.8]{eckl}, we have
$$
 n(L) \ge \nu (L) \ge \kappa (L),
$$
 where $\kappa (L)$ is the Kodaira dimension of $L$
 and $\nu (L)$ is the numerical Kodaira dimension of $L$.
 Thus the assumption $n(L) < \dim X$ implies $\nu (L) < \dim X$.
 We note that if $X$ is an rreducible symplectic manifold
 and $L$ a nef divisor on $X$,
 then $q(L) = 0$ if and only if $\nu (L) < \dim X$ by 
 \cite[Theorem 4.7]{fujiki}.
\end{remark}

\begin{remark}
 If we consider an irreducible symplectic manifold which
 is the moduli space of semi stable torsion free sheaves on 
 a $K3$ surface or an abelian surface, 
 then we have other existence
 conditions of fibre space structures.
 Please see \cite[Theorem 1.3]{gulbranden}, \cite[Theorem 4.3]{markshevich} and
 \cite[Theorem 2]{sawon}.
 On the other hand, Amerik and Campana
 study a rational map $f : X \dasharrow S$ from an
 irreducible symplectic manifold $X$ to a normal variety
 $S$ such that $0 < \dim S < \dim X$ and the Kodaira
 dimension of a general fibre of $f$ is zero.
 Please see \cite[Th\'eor\`eme 3.6]{ekaterina}.
\end{remark}

\noindent
{\bfseries Acknowledgements. \quad}
 The author express his thanks to Professor O.~Fujino who
 informed me the nef dimension. He also express his thanks 
 to Professor E.~Amerik who sended me a copy of her preprint.
 By that article, the proof of Theorem \ref{main} is simplified.
 He also express his thanks to Professor S.~Boucksom who
 gave me an important comment of Claim \ref{key}.

\section{Proof of Theorem \ref{main}}

\begin{say}
 By \cite[Theorem 1]{kawamata}, it is enough
 to prove that
$$
 \nu (L) = \kappa (L).
$$
 By the assumption, there exists a rational map
 $f : X \dasharrow S$ which satisfies 
 the three conditions of Definition--Theorem \ref{Nef_Dimension}.
 Let $\pi : Y \to X$ be a resolution of indeterminacy of $f$.
 We denote by $g$ the induced morphism. 
$$
 \xymatrix{
   X  \ar@{-->}[d] & Y \ar[l]_{\pi} \ar[ld]^g \\
       S               &
 }
$$
\end{say}
\begin{claim}\label{Claim1}
 Let $L_S$ be a very ample divisor on $S$. Then
$$
 q(\pi_* g^* L_S ) = 0.
$$
\begin{proof}
 Let $H_1$ and $H_2$ be general members of
 the linear system of $g^* L_S$.
 By \cite[Theorem 5]{beauville},
$$
 q (\pi_* g^* L_S) = \frac{\dim X}{2}
 \int_{X} \pi_* H_1
 \pi_* H_2
 (\omega \bar{\omega})^{(1/2)\dim X - 1} ,
$$
 where $\omega$ is a symplectic form on $X$.
 Since $L_S$ is very ample, 
 the intersection $\pi_* H_1$ and $\pi_* H_2$
 defines a codimension $2$ effective cycle.
 Thus $q (\pi_* g^* L_S) \ge 0$.
 We derive a contradiction assuming that
 $q (\pi_* g^* L_S) > 0$. 
 Let $\mathcal{P}$ be the positive cone attached to
 $q$ and $\mathcal{E}$ the pseudo effective cone.
 By \cite[Theorem 4.3 (i)]{boucksom},
 $\mathcal{P} \subset \mathcal{E}$.
 Hence  $\pi_* g^* L_S$ defines an interior point
 of $\mathcal{E}$. This implies $\pi_* g^* L_S$ is big.
 Let $F$ be a general fibre of $g$.
 Then $\pi^* \pi_* g^* L_S|_{F}$ is big.
 Since $\pi$ is isomorphic in some neighbourhood of  $F$,
 $\pi^* \pi_* g^* L_S |_{F} = g^* L_S |_{F} = 0$.
 That is a contradiction.
\end{proof}
\end{claim}
\begin{claim}\label{Claim2}
$$
 q(L, \pi_* g^* L_S) = 0.
$$
\end{claim}
\begin{proof}
 Since $L$ is nef, $q (L, \pi_* g^* L_S) \ge 0$.
 We compute $q (L + \pi_* g^* L_S)$.
 By the assumption and Remark\ref{relation_of_dimensions},
 $q(L) = 0$. By Claim \ref{Claim1},
 $q(\pi_* g^* L_S) = 0$.
 Thus
$$
 q (L + \pi_* g^* L_S) = 2 q (L, \pi_* g^* L_S ).
$$
 If we assume  $q (L, \pi_* g^* L_S) > 0$,
 then we have $L + \pi^* g_* L_S$ is big
 by the same argument in the proof of Claim \ref{Claim1}. 
 However
$$
 \pi^*(L + \pi_* g^* L_S) |_{F} = \pi^* L \equiv 0
$$
 by the condition (2) of Definition--Theorem \ref{Nef_Dimension}.
 That derives a contradiction.
\end{proof}

\begin{claim}\label{key}
 There exists a positive constant $\lambda_0$ such that
$$
 \lambda_0 L \sim_{\mathbb{Q}}  \pi_* g^* L_S.
$$
\end{claim}
\begin{proof}
 Let $L^{\vee}$ be the hyperplain in 
 $H^{1,1}(X,\mathbb{C})_{\mathbb{R}}$
 defined by $q(L,*)$.
 Assume that $L^{\vee}$ contains an interior of the
 positive cone $\mathcal{P}$ of $X$.
 Then there exists a divisor $H$ such that
 $q(H) > 0$ and $q(L,H) < 0$.
 However this derives a contradiction because
 the signature of $q$ on  $H^{1,1}(X,\mathbb{C})_{\mathbb{R}}$
 is $(1,k - 1)$, where $k = \dim H^{1,1}(X,\mathbb{C})_{\mathbb{R}}$.
 Hence $L^{\vee}$ intersects only the boundary of
 the closure $\bar{\mathcal{P}}$ of $\mathcal{P}$.
 Since $\mathcal{P}$ is a quadratic cone,
 $L^{\vee} \cap \bar{\mathcal{P}}$ coincides with the ray which is
 generated
 by $L$. By Claim \ref{Claim2},
 $L$ and $g^*L_S$ are contained in the same ray.
\end{proof}
\begin{say}
\begin{proof}[Proof of Theorem \ref{main}.]
 By Claim \ref{key}, we obtain 
$$
 \lambda_0 \pi^* L \sim_{\mathbb{Q}} g^* L_S + \sum e_i E_i
$$
 where $E_i$ is a $\pi$-exceptional divisor.
 Hence we have the following inequalities:
$$
  \nu (L) \ge \kappa (L) \ge \kappa (g^* L_S) = n (L).
$$
 By the inequality in Remark \ref{relation_of_dimensions}, we are done.
\end{proof}  
\end{say}

\end{document}